\documentclass[12pt]{amsart}
\usepackage{amssymb,amsmath}
\usepackage{epsf,latexsym,graphicx}
\usepackage{psfrag}
\psfragscanon
\theoremstyle{plain}
\newtheorem{thm}{Theorem}[section]

\newtheorem{cor}[thm]{Corollary}
\newtheorem{lemma}[thm]{Lemma}
\newtheorem{example}{Example}
\theoremstyle{definition}
\newtheorem{definition}{Definition}
\newtheorem{remark}{Remark}
\numberwithin{equation}{section}

\begin{document}

\title{Reflections on symmetric polynomials and arithmetic functions}
\author{Trueman MacHenry}
\author{Geanina Tudose}
\begin{abstract}
In an isomorphic copy of the ring of symmetric polynomials we study some families of
  polynomials which are indexed by rational weight vectors. These
  families include well known symmetric polynomials,
  such as the elementary, homogeneous, and power sum symmetric
  polynomials. We investigate properties of these families and focus on
   constructing their rational roots under a
  product induced by convolution. A direct application of the latter
  is to the
  description of the roots of certain  multiplicative arithmetic
  functions (the core functions) under the convolution product.
\end{abstract}

\address{Department of Mathematics and Statistics\\
York University\\
North York, Ont., M3J 1P3\\
CANADA}
\email{Trueman.Machenry@mathstat.yorku.ca}
\email{gtudose@mathstat.yorku.ca}

\maketitle

\section{Introduction}
 
This paper is concerned with a certain isomorphic copy of the ring 
$\Lambda \otimes_{\mathbb{Z}} R$ of 
symmetric functions, namely the ring of isobaric polynomials 
$\Lambda'$, where the 
isomorphism is given by a  polynomial map involving the 
elementary symmetric functions. The ring will be taken to be 
either the integers $\mathbb{Z}$ or the rationals $\mathbb{Q}$, and the 
image of a symmetric polynomial under the isomorphism mentioned above will 
be called an isobaric reflect. An isobaric\footnote[3]{the term  isobaric
  is due to P\'{o}lya~\cite{P}; the cycle index of a finite group
  appearing in P\'{o}lya's Counting Theorem is an isobaric polynomial.}
polynomial is one of the form 
$P_n= \sum_\alpha A(\alpha) t_1^{\alpha_1}...t_k^{\alpha_k}$, where  
$\alpha= (\alpha_1,...\alpha_k),\  \alpha_i \geq 0$ are integers with $\sum_j 
j\alpha_j =n$. Thus such a polynomial can be represented as a polynomial 
in the Young diagrams of $n$, where the monomial is an encoding of the partition of 
$n$. As for the ring of symmetric polynomials, we can allow either a finite 
number $k$ of variables or we can work in $\oplus_k \Lambda_k'$ with 
infinitely many variables.

Families of isobaric polynomials occur in many contexts in mathematics. 
In~\cite{TM1} it was shown that the reflects of the complete symmetric polynomials (CSP) 
determine the multiplicative arithmetic functions locally. In~\cite{TM2} 
it was shown 
that the reflects of the power sum  symmetric polynomials (PSP) determine the 
lattice of root fields of quadratic extensions. Properties of these two 
sequences of polynomials were discussed in~\cite{TM3} where the 
CSP-reflects are 
called Generalized Fibbonacci Polynomials (GFP), and the PSP-reflects are 
called the Generalized Lucas Polynomials (GLP). Recall that the Complete 
Symmetric Polynomials form a $\mathbb{Z}$-basis for the symmetric 
polynomials as do the Elementary Symmetric Polynomials (ESP), while the 
Power Symmetric Polynomials (PSPs) form a $\mathbb{Q}$-basis. The 
analogues of these facts carries over to the isobaric polynomial algebras 
by way of a canonical isomorphism of the ring $\Lambda$ to the ring 
$\Lambda'$ denoted by $\varXi$. In fact, this 
isomorphism is just the one that takes symmetric functions on $k$ 
variables, written in terms of elementary symmetric polynomials, and 
rewrites each elementary polynomial $e_j$ as $(-1)^{j+1}t_j$.
$$
 \hat{e}_i= \varXi(e_i)= (-1)^{j+1}t_j.
$$

It is well-known that the Schur Symmetric Polynomials (SSP) determine the 
complex character table of the finite symmetric groups using the 
Littlewood-Richardson rule and the Frobenius Character Theorem. Thus the 
SSPs for a given $n$ can be regarded as an encoding of the complex 
character table of $Sym(n)$. The Frobenius Character Theorem can be 
written in terms of isobaric polynomials, namely in terms of GLPs. Using 
this fact, the complex characters of $Sym(n)$ can be easily calculated 
from the isobaric reflects of the SSP$(n)$, the Schur polynomials, for a 
given $n$.

The families GFP and GLP have the additional useful property that each 
satisfies recursion relations (Newton identities). It will turn out that  
$\Lambda'$ contains a large class of recursively defined families 
(Theorem~\ref{free}). These are the families of what we have called {\it 
weighted 
isobaric polynomials} (WIPs), and they are the main subject of this paper. 
Such polynomials are determined by assigning a weight to each of the 
variables $t_j$, i.e. by assigning a weight vector to the set of variables 
$\{ t_j \}$. Such families will be called weighted isobaric families. It 
turns out that the union of all such families does not exhaust the ring of 
isobaric polynomials. In fact, in this paper we show that among the Schur 
polynomials, exactly those Schur reflects which represent hook (Young) 
diagrams can belong to a sequence of weighted isobaric polynomials 
(Theorem~\ref{rechooks}  and Theorem~\ref{rechooks2}).

Families of WIPs, multi-indexed by their weights, form in a natural way a 
free abelian group induced by addition of their weight vectors 
(Theorem~\ref{free}). 
The weighted families GFP and GLP, i.e. the CSP and PSP reflects, are the 
weighted families determined by weight vectors $(1,1,...)$ in the case of 
GFP, and $(1,2,3,...)$ in the case of GLP. The Schur-hook reflects have 
weight vectors of the form $(0,0...1,1,...)$ (Theorem~\ref{hook}).

The coefficients of the monomials in an isobaric polynomial are uniquely 
determined by the exponents of the variables and the weight vector of the 
family (Theorem~\ref{coeffthm}). In order to prove Theorem~\ref{coeffthm} 
we use the fact 
that each monomial 
determines a lattice whose nodes are the Young diagrams of the monomial 
obtained by derivation. However, this lattice is not the well-known 
Young's 
lattice, but instead it is a lattice partially ordered by the
pointwise inequality of the exponents of the constituent nodes. It assumes a 
major role in this paper in understanding the construction of the WIPs and 
certain other structures associated to them. 

In 1988, Carrol and Gioia~\cite{CG} gave a numerical description of the $q$-th roots ($q \in \mathbb{Q}$) of the multiplicative arithmetical functions in the 
group of units of the ring of arithmetical functions. In~\cite{TM1} it was 
shown that 
under  convolution, these functions form a free abelian group 
generated by the completely multiplicative functions, as mentioned above; 
it was also shown in that paper that the GFPs give a generic set of 
generating functions for this group of arithmetic functions in the 
following sense: each multiplicative function in the core
group\footnote[4]{The {\sl Core} group is the subgroup of the group of units in
  the ring of arithmetic functions generated by the complete
  arithmetic functions.} of the 
group of units of a multiplicative arithmetic function together with its 
convolution inverse is uniquely determined locally by a monic polynomial (over the 
complex fields), the generating polynomial. What is called a negative 
element is a multiplicative function whose local values are just the coefficients 
of this generating polynomial, while the inverse of this negative core 
function, a positive element, is a multiplicative function whose local
values 
are 
given by evaluating the series of GFPs truncated at the degree of the 
generating  polynomial at these coefficients. In Section 4, we produce a 
sequence of isobaric polynomials which are the $q$-th roots for any $q 
\in \mathbb{Q}$ of 
the generic generating functions for these roots, that is, the $q$-th 
roots of the Generalized Fibonacci and Lucas polynomials (the CSP and PSP 
reflects). Thus, we have embedded the core group into its divisible 
closure.

Moreover, our construction is more far-reaching than this. It produces a 
set of $q$-th roots with respect to a product induced by convolution, which we have called the 
level product (so-called because it acts on polynomials of the same level 
and conserves level) for every isobaric 
polynomial in any weighted family (Theorem~\ref{hsum}). Theorem~\ref{hsum} 
also implies that 
under the level product, an element has a level product inverse.

It is not the case that the isobaric roots of weighted functions are 
necessarily weighted.  However, they are determined by specifying a weight 
vector. The appropriate structure to look at here, then, is the ring 
generated by the level product. This is a graded ring $\mathcal{H}$ 
containing WIPs. However, there is still more algebraic structure. Since we have inversion under 
level-product, each 
weighted family and its divisible closure has the structure of the rationals. Moreover, because of 
Theorem~\ref{free} each weighted family is acted on by translations, and 
so we 
finally have that all of this together with the derivation operation give 
a differential graded abelian group acted on by an affine group.
 
\section{Weighted Isobaric Polynomials}

\begin{definition}
A {\it weighted isobaric polynomial (WIP)} is defined inductively in terms
of a particular lattice which will now be described. The nodes of the
lattice are the monomials $\{ t_1^{\alpha_1}\ldots t_k^{\alpha_k}\}$,
where $\alpha=(\alpha_1,\ldots \alpha_k),\  \alpha_j \in \mathbb{Z},\ \alpha_j \geq 0$ 
with $\sum_j j\alpha_j=n$, for a fixed $n$. We denote this by $\alpha 
\vdash n$ and by $|\alpha|$ the sum $\sum_j \alpha_j$.
The relation
$$
 t^\beta=(t_1^{\beta_1}\ldots t_k^{\beta_k}) \leq
 t^\alpha=(t_1^{\alpha_1}\ldots t_k^{\alpha_k}),\ \mbox{if } \beta_j
 \leq \alpha_j, \ \ \mbox{for every }j
$$
\noindent
imposes a lattice structure on the set of $\{t^\alpha\}$. The {\it
  depth} of a monomial in the lattice is $(\sum \alpha_j)$, and its
  {\it level} is $n=\sum_j j\alpha_j$. Let $1$ be the bottom element. We
  assign a weight vector $\omega= (\omega_1,\ldots \omega_k)$ to the
  variables $(t_1 \ldots t_k)$, i.e. to each $t_j$ we assign
  $\omega_j$ as its coefficient. The weights are taken to be rational
  but they can belong to any ring $R$. Having done this, $t^\alpha$ will be
  assigned the coefficient equal to the sum of the coefficients  of
  all $t^\beta$ that have depth $(\sum \alpha_j-1)$ and for which 
$t^\beta < t^\alpha$. Thus each monomial involving $t^\alpha$ can be
  associated with a (finite) sublattice $\mathcal{L}(t^\alpha)$ with
  top element $t^\alpha$, the lattice of all those monomials whose
  coefficients contribute to the coefficient of $t^\alpha$. For any
  two monomials there is a monomial whose lattice contains their
  lattices\footnote[5]{This lattice can be thought as a lattice of Young
  diagrams $(1^{\alpha_1},..k^{\alpha_k})$ in which a ``smaller''
  diagram is one with one less row; it is clearly not a Young
  lattice. As far as we know these lattices have not been introduced
  before in the study of symmetric functions.}. Note also that the lengths of all maximal chains of the sublattice
$\mathcal{L}(t^\beta)$ with $t^\beta < t^\alpha$ 
are the same and clearly equal to the corresponding {\it depth}.

 The isobaric polynomial of level $n$ and weight $\omega$
  is therefore $P_{k,n,\omega}=\sum_{\alpha \vdash
  n}weight(t^\alpha)t^\alpha$. 
\end{definition}

\vskip12pt
\begin{example}
The lattice $\mathcal{L}(t_1^2t_2t_3)$ is 
\begin{figure}[h]
\includegraphics{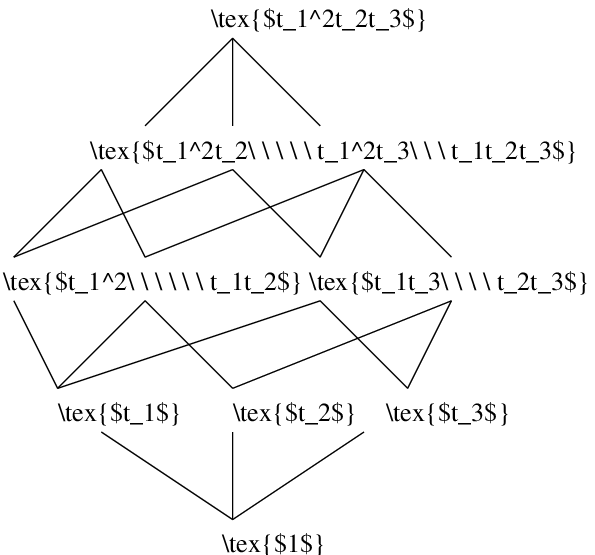}
\end{figure}
\end{example}

We can now state

\begin{thm}\label{coeffthm}
If $\omega$ is a weight vector, then the WIP of degree $n$ has at most
$\mathcal{P}(n)$ terms, where $\mathcal{P}(n)$ is the number of partitions of $n$, and
  the coefficients are given by 
\begin{equation}\label{coeff}
weight(t^\alpha) = A_{k,n,\omega}(\alpha) ={{\sum \alpha_i} \choose {\alpha_1,...\alpha_k}}\frac{\sum_i\alpha_i\omega_i}{\sum_i \alpha_i}
\end{equation}
\end{thm}

\noindent
In particular, the coefficients of the families $\mathcal{F}_k$, the GFP and
$\mathcal{G}_k$, the GLP, are given by $\displaystyle{ {{\sum \alpha_i} \choose
  {\alpha_1,...\alpha_k}} }$ and 
$\displaystyle{
  n\frac{[(\sum_{j=1}^{k}\alpha_j-1)]!}{\prod_{j=1}^k(\alpha_j)!}}$,
respectively, where the weight vector of the GFPs is given by
$(1,1,...)$ and the weight vector of the GLPs is given by
$(1,2,3,...)$. It is of interest that, when $\omega$ is an integer
vector, the numbers $A_{k,n,\omega}(\alpha)$ are integers. This
 will be a trivial consequence of the proof of Theorem~\ref{coeffthm}.

\noindent
{\bf Proof of Theorem~\ref{coeffthm}}.

Let $\omega =(\omega_1,\ldots \omega_k)$ with $\omega_j \in \mathbb{Z}$ be a weight assignment
to the indeterminates $t_1, \ldots t_k$. This assignment, together
with the inductive rule for determining the coefficient of a monomial
in the lattice, will define a family of WIPs, denoted by
$\mathfrak{F}_{k,\omega}$ or just $\mathfrak{F}_{\omega}$.

To see that the coefficients are as stated in the theorem, we proceed
by induction on the depth to compute the coefficient of $t^\gamma$, 
where $t^\gamma= t_1^{\gamma_1}...t_k^{\gamma_k}$. The monomials that
contribute to the coefficient of $t^\gamma$ are just
$t^{\gamma^{(j)}}$, where $\gamma^{(j)}= (\gamma_1, \ldots
\gamma_j-1,\ldots\gamma_k)$. Then by induction
$$
A_{k,n,\omega}(\gamma^{(j)})  = \frac{[(\sum_i \gamma_i)-2]!}{\prod_{i \neq j}(\gamma_i)!(\gamma_j-1)!}[\sum_{i \neq j}\gamma_i\omega_i+(\gamma_j-1)\omega_j]
$$

\noindent
and so 
$$
A_{k,n,\omega}(\gamma) = \sum_{j=1}^k \frac{[(\sum_i\gamma_i)-2]!}{\prod_{i\neq j}
(\gamma_i)!(\gamma_j-1)!}[\sum_{i \neq j}\gamma_i\omega_i+(\gamma_j-1)\omega_j]=
$$
$$
= \sum_{j=1}^k \frac{[(\sum_i\gamma_i)-1]!(\gamma_j)![\sum_i
  \gamma_i\omega_i-\omega_j]}{\prod_{i}(\gamma_i)! \sum_i
  (\gamma_i-1)(\gamma_j-1)!} =
\frac{((\sum_i\gamma_i)-1)!}{\prod_i(\gamma_i)!}
\Big{(}\frac{\sum_i\gamma_i\omega_i}{(\sum_i\gamma_i)-1}\sum_{j=1}^k
  \gamma_j- \frac{\sum_i\gamma_i\omega_i}{(\sum_i\gamma-1)}\Big{)}
$$
$$
=\frac{((\sum_i\gamma_i)-1)!}{\prod_i(\gamma_i)!}
(\sum_i\gamma_i\omega_i).
$$
\hfill
$\Box$

\noindent
Note that a family $\mathfrak{F}_{k,\omega}$ is a {\it sequence} of
  polynomials, one for each degree. For example the first four of
  these in any sequence are of the form:

$$
\begin{array}{lll}
P_{1,\omega} & = & \omega_1t_1\\
P_{2,\omega} & = & \omega_1t_1^2 +\omega_2t_2\\
P_{3,\omega} & = & \omega_1t_1^3+ (\omega_1+\omega_2)t_1t_2
+\omega_3t_3\\
P_{4,\omega} & = & \omega_1t_1^4 +
(2\omega_1+\omega_2)t_1^2t_2+\omega_2t_2^2+(\omega_1+\omega_3)t_1t_3+\omega_4t_4
\end{array}.
$$
\noindent
Since all of the operations involved in computing the coefficients are
ring operations, we have the following.

\begin{cor}
If $\alpha$ and $\omega$ are integer vectors, then $A_{k,n,\omega}(\alpha)$ is an integer.
\end{cor}

\begin{remark}
For a weight $\omega$ with $\omega_i \neq 0$, for any $i$,  using the jacobian criterion we have that 
$Jacob(P_{i,\omega})=\prod_i\omega_i \neq 0$ and thus 
the family $\mathcal{F}_{k,\omega}$ consists of algebraically independent 
polynomials. This allow us to construct a new basis for $\Lambda'$ for 
each such weight vector $\omega$, whose elements are
$$
    P_{\lambda,\omega}= \prod_j P_{\lambda_j,\omega}, \qquad \mbox{for 
every partition } \lambda.
$$
\end{remark}

\noindent
Moreover, if $P_{k,\omega} \in \mathfrak{F}_{k,\omega}$ and,
$P_{k,\omega'} \in \mathfrak{F}_{k,\omega'}$ then $P_{k,\omega} + P_{k,\omega'} \in \mathfrak{F}_{k,\omega+\omega'}$. If we define addition on
these classes by
$\mathfrak{F}_{k,\omega}+
\mathfrak{F}_{k,\omega'}=\mathfrak{F}_{k,\omega+\omega'}$ we have
\vskip12pt
\begin{thm}\label{free}
$\{ \mathfrak{F}_{k,\omega} \}$ is a free $\mathbb{Z}$-module under
  this operation.
\end{thm}

\noindent
{\bf Proof}
First notice that $A_{k,n,\omega}(\alpha)=A_{k,n,\omega'}(\alpha)$ if and only if $\omega= \omega'$, for by
Theorem~\ref{coeffthm} this implies that $\sum_j \alpha_j\omega_j=
\sum_j \alpha_j\omega_j'$, and this in turn implies, by taking
different values for the $\alpha$'s, that $\omega_j=\omega_j'$. Thus
equality of families implies identity of weights. Clearly the set of 
families is a
$\mathbb{Z}$-module under the defined operation. Since there is a homomorphic
image of the $\mathbb{Z}$-module of families onto the $\mathbb{Z}$-module of integer weight
vectors, the assertion follows.

\hfill
$\Box$
\vskip12pt
While the WIPs form a graded (and, as we shall see, a differential graded)
ring, they are not closed under multiplication, so the best we can
do is speak of the subring of $\Lambda'$ generated by WIPs. We shall
see later that this is a proper subring of $\Lambda'$. The lattice
$\mathcal{L}(t^\alpha)$ is generated by derivations $\displaystyle{
  \frac{1}{\alpha_i}\partial_i}$, where $\partial_i:t^\alpha \rightarrow
  \alpha_i t^\beta$, and where $\beta= (\alpha_1,...\alpha_i-1,...\alpha_k)$. 
Thus $\Lambda'$ becomes a differential graded ring. We
    also shall need the {\it total} differential operator $D_j$ where
    $D_j=D_1(D_{j-1})$ and $D_1(t^\alpha)= \sum_i \partial_i(t^\alpha)$. We
    are interested in the total differential when it is evaluated at a
    weight vector $\omega=(\omega_i)_i$. We write this as
    $D_j(\omega^\alpha)$. We shall need the following lemma later.

\begin{lemma}\label{ll}
$$ D_{(\sum \alpha_i)}(\omega_1^{\alpha_1}\ldots \omega_k^{\alpha_k})=
(\sum_{i=1}^k \alpha_i-1)!(\alpha_1\omega_1+\ldots +\alpha_k\omega_k).
$$
\end{lemma}

(Note that the right hand side of the Equation above is just
$\prod(\alpha_i)!$ times the coefficient of a monomial term in a
weighted isobaric polynomial given in Theorem~\ref{coeffthm}.)

\noindent
{\bf Proof of Lemma~\ref{ll} } 
Let $u= \sum\alpha_i-1$, then
$$
D_u(\omega_1^{\alpha_1}\ldots \omega_k^{\alpha_k}) = \sum_j \alpha_j D_{u-1}(\omega_1^{\alpha_1}\ldots
\omega_j^{\alpha_j-1}\ldots \omega_k^{\alpha_k} )
$$
\noindent
and by induction it is 
$$
= (\sum_i\alpha_i-2)!(\sum_j\alpha_j(\sum_i((\alpha_i\omega_i)-\omega_j)))= 
(\sum_i\alpha_i-2)!(\sum_i(\alpha_i\omega_i)(\sum_j\alpha_j)-(\sum_i\alpha_i))
$$
$$
=(\sum_i\alpha_i-2)!(\sum_i \alpha_i-1)(\sum_i \alpha_i\omega_i)=
(\sum_i \alpha_i-1)!(\sum_i \alpha_i\omega_i).
$$
\hfill
$\Box$

\section{Weighted isobaric polynomials and Schur  functions}

Denoting the  Schur polynomials by SSP, we want to consider the  family of isobaric reflects
of the SSP, the Schur reflects. In this section we are interested in the
question of which Schur reflects can belong to a sequence of
WIPs. More generally, of course, we want to know when any isobaric polynomial
belongs to some WIP sequence. We have solved the problem for Schur
reflects, but we still have no interesting criterion in the general case.

First, a look at the classical role that Schur polynomials play from
the point of view of reflects. We have two routes from the ring of symmetric polynomials $\Lambda$ to the 
ring
of isobaric polynomials $\Lambda'$, one ``backward'' from the character
table of $Sym(n)$ by way of the Frobenius Character Theorem; the other
``forward'', using the Jacobi-Trudi formulae. Recall that the
Frobenius theorem for $S_\lambda$, where $\lambda$ represents a
conjugacy class of $Sym(n)$ (or equivalently a Young diagram) is the
formula
$ S_\lambda= \frac{1}{n!}\sum_\mu C(\mu)\chi_\lambda^\mu P_\mu$, where
$P_\mu$ is the power symmetric polynomial, $C(\mu)$ the size of the
conjugacy class of $\mu$, $S_\lambda$ the Schur polynomial of shape
$\lambda$  and $\chi_\lambda^\mu$ is the character of $Sym(n)$
afforded by $\lambda$ applied to $\mu$. 

Using the bijective mappings $\varXi_k$ which takes $G_{k,n}$'s to
PSPs, i.e to $P_\mu$'s, we have an analogue of the Frobenius Character Theorem in terms
of reflects relating the character afforded by $\lambda$ to the PSP
isobaric polynomials, namely $ \hat{S}_\lambda= \frac{1}{n!}\sum_\mu
C(\mu)\chi_\lambda^\mu G_\mu$. This gives us a representation of the
complex character table of $Sym(n)$ in terms of PSP reflects, the $G_\mu$'s. Rewriting these polynomials in terms of the CSPs using,
e.g., Theorem 3 of~\cite{TM1} gives a representation in terms of
Complete Symmetric Polynomials.

The ``forward'' route using the Jacobi-Trudi
identities, writes the Schur reflects either in terms of the CSP or in terms of the
ESP, that is in terms of the variables $t_j$. This route gives us
determinantal formulae for the Schur reflects either in terms of the
polynomials $F_{k,n}$ or in terms of the basic variables. Thus
$\hat{S}_\lambda = det[F_{\lambda_j-i+j}]_{1\leq i,j\leq k} =
det[t_{\lambda_j-i+j}]_{1\leq i,j\leq k}$. These maps reveal an
interesting iterative property of those Schur polynomials determined
by partitions $\lambda$ which are hooks, i.e. of type $\lambda=(p,
1^q)$. The hook Schur polynomials can be expressed in terms of the
homogeneous and elementary symmetric polynomials as
$S_{(n-r,1^r)}=\sum_{j \geq r+1}(-1)^{j-r-1} e_jh_{(n-j)}$ [\cite{M}
Chap 3]. Taking the isobaric reflect we obtain 
\begin{thm}\label{rechooks}
$$
\hat{S}_{(n-r,1^r)}= (-1)^r \sum_{j=r+1}^{n}t_j\hat{S}_{(n-j)}, \qquad
0\leq r \leq n
$$
\end{thm}

\noindent
Now it is clear that $\hat{S}_{(n)}=F_n$, hence Theorem~\ref{rechooks}
can be written as
\begin{thm}\label{rechooks2}
\begin{equation}\label{sumhooks}
\hat{S}_{(n-r,1^r)}= (-1)^{r} \sum_{j \geq r+1}^{n}t_j F_{(n-j)}.
\end{equation}
\end{thm}

The reflects of the Schur polynomials determined by hooks form families of 
weighted isobaric polynomials. If we denote the horizontal boxes of the 
hook as the {\it arm} of the hook and the vertical boxes as the {\it leg} of the 
hook then those Schur polynomial reflects determined by hooks with legs of 
the same length belong to the same weighted family, where the length is 
the number of boxes. More precisely we have

\begin{thm}\label{hook}
Hooks with leg length $(r+1)$ belong to the weighted family determined by the weights $\omega_{(r)}=(-1)^{r}(0,...0,1,1,...)$, with $r\ 0$'s, the 
rest $1$'s.
\end{thm}

\noindent
{\bf Proof of Theorem~\ref{hook} }
We will first show that each term on the right hand side
of~(\ref{sumhooks})  is a
weighted polynomial and describe its weight.

\begin{lemma}
The weight vector for $t_jF_{n-j}$ is $(0,...0,1,0,...)$, where the $j$-th component is $1$.
\end{lemma}

\noindent
{\bf Proof}
Exponents of the monomials in $F_{n-j}$ satisfy $\sum i \alpha_i= n-j$. So 
the coefficients of $F_{n-j}$ are by Theorem~\ref{coeffthm}, $A_{k, n-j,(1,1,..)}(\alpha)=\displaystyle{ 
\frac{(\sum_{i=1}^k\alpha_i)!}{\prod_{i=1}^k (\alpha_i)!} }$, where $\sum i \alpha_i= n-j$.
But the coefficient of the monomial in $t_jF_{n-j}$ with exponent $\beta= 
(\alpha_1,...\alpha_j+1,... \alpha_k)$, that is the coefficients in 
$A_{k,n,\omega}(\beta)$, is the same, i.e  $A_{k,
  n-j,(1,1,...)}(\alpha)=A_{k,n,\omega}(\beta)$. On the other hand
$$
A_{k,n,\omega}(\beta)=
\frac{[(\sum_{i=1}^k\alpha_i)]!}{(\alpha_j+1)\prod_{i\neq
    j}(\alpha_i)!}[\sum_{i\neq j} \alpha_i\omega_i+(\alpha_j+1)\omega_j].
$$
From these considerations we get that $\sum_{i\neq j}
\alpha_i\omega_i+(\alpha_j+1)\omega_j=\alpha_j+1$. Since this is true for all 
exponents such that $\sum i \alpha_i= n-j$ we must have that the
weight vector of $t_jF_{n-j}$ is $(0,...0,1,0,...)$.
\hfill
$\Box$

We complete the proof of Theorem~\ref{hook} by using the result of
Theorem~\ref{free} which shows that a sum of weighted polynomials is a
weighted polynomial with weight the sum of weights of each term.
Therefore the proof is complete. 
\hfill
$\Box$
\vskip12pt

Next we give a complete answer to the question of which Schur polynomials 
belong to  weighted families.

\begin{thm}\label{nonhooks}
If $\lambda$ is a partition of $n$ such that the shape of $\lambda$ is not 
that of a hook, then $\hat{S}_\lambda \neq P_{n,\omega}$, for any weight 
vector $\omega$. The Schur reflect cannot be a weighted polynomial.
\end{thm}

Before we proceed with the proof we need

\begin{definition}[lexicographic order on $\mathcal{P}(n)$]
 
Let $\lambda=(\lambda_1\geq \lambda_2 \geq ...)$ and $\mu=(\mu_1 \geq
\mu_2 \geq ...)$ be two partitions of $n$. We say that $\lambda
< \mu$ if and only if, for some index $i$ we have $\lambda_j = \mu_j$
for $j< i$ and $\lambda_i <\mu_i$. 
\end{definition}

\begin{example}
The lexicographic order on $\mathcal{P}(4)$.
$$
(1^4) < (1^2,2) < (2^2) < (1^1, 3^1) < (4^1)
$$
\end{example}

This order induces a corresponding order on the monomials $t^\alpha$ with 
$(1^{\alpha_1},...,k^{\alpha_k}) \vdash n$. Furthermore we shall write 
the WIPs by ordering its monomials starting with the smallest. For example $P_{4,\omega}= \omega_1t_1^4 +
(2\omega_1+\omega_2)t_1^2t_2+\omega_3t_2^2+(\omega_1+\omega_3)t_1t_3+\omega_4t_4$.

Next we place in a \fbox{{\it box $i$}} all monomials $t^\alpha= 
t_1^{\alpha_1}\ldots t_i^{\alpha_i}$ such that $\alpha_i \neq 0$ and
$\alpha_j =0$, for $j >i$. 
\begin{example} The arrangement of boxes in $\mathcal{P}(4)$.
\end{example}
\begin{figure}[h]
\includegraphics{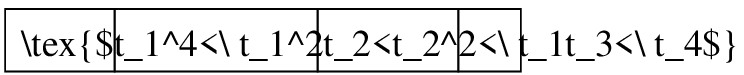}
\end{figure}
\noindent
It is easy to see  that ordering the boxes according to their indices,  
gives a saturated chain under the lexicographic order of all  monomials 
$t^\alpha$, with $(1^{\alpha_1},...,k^{\alpha_k}) \vdash n$.
We note that the smallest monomial in each box corresponds to a hook
and only to a hook. This follows from the way we defined the boxes: the
smallest monomial in, say \fbox{box $i$} is $t_1^{n-i}t_i$ which
corresponds to the hook $(i,1^{n-i})$. 

\noindent
{\bf Proof of Theorem~\ref{nonhooks}}
Via the Jacobi-Trudi identity we have $S_\lambda= det(e_{\lambda_i'-i+j})$ 
where $\lambda'$ is the conjugate partition of $\lambda$ obtained by 
transposing the Young diagram.
Under the reflection isomorphism, i.e. $\hat{e}_i= (-1)^{i-1}t_i$, we 
get $\hat{S}_\lambda= det((-1)^{\lambda_i'-i+j-1}t_{\lambda_i'-i+j})$ .
In the expression above the smallest monomial is obtained from the main 
diagonal of the determinant (as the transition matrix from the bases
$S_\lambda$ and $e_\lambda$ is upper triangular [Chap 6 in~\cite{M}]
and also the Appendix). 
The smallest monomial is $t^\delta$ where 
$(1^{\delta_1},...s^{\delta_s})=\lambda'$ and its coefficient is 
$(-1)^{n-\lambda_1} \neq 0$. 
Assume now that there exists a weight $\omega$ such that
$\hat{S}_\lambda = P_{n,\omega}$. Recall that from 
Theorem~\ref{coeffthm}, in $P_{n,\omega}$
the coefficient of $t^\alpha$ is $A(\alpha)= {{\sum \alpha_i} \choose {\alpha_1,...\alpha_k}}\frac{\sum_i\alpha_i\omega_i}{\sum_i \alpha_i}$.
Assume that the smallest monomial $t^\delta$ belongs to \fbox{box
  $s$}, i.e, 
$\delta_s\neq 0$ and $\delta_i=0$ for $i>s$. Since $\lambda$ is not a 
hook, $\lambda'$ is not a hook either and so $t^\delta$ is not the first monomial in box $s$. Moreover, we must have that $A(\beta)=0$ for any 
$\beta$ such that $(1^{\beta_1},...k^{\beta_k}) <
(1^{\delta_1},...s^{\delta_s})$  in the lexicographic order. In particular 
$A(\beta^{(i)})=0$, for $\beta^{(i)}=(n-i,0,...1...0)$, with $1$ in the 
$i$-th place, for $i=1,2,... s$,  that is, the first monomials in boxes 
$1,2,...s$. This is to say
$$
  {{n-i+1}\choose {1}} \frac{(n-i)\omega_1 + 1\omega_i}{n-i+1} = 0, \quad i=1,\ldots s.
$$

\noindent
From this system of equation we get 
$$
 \omega_1=\omega_2=\ldots=\omega_s=0
$$
which in turn gives $A(\delta) ={{\sum \delta_i}\choose 
{\delta_1,\delta_2,...}}\frac{ \sum_{i=1}^s \delta_i\omega_i}{\sum \delta_i}=0$, a contradiction.
\hfill
$\Box$
\vskip12pt
Theorem~\ref{hook} and~\ref{nonhooks} and tell us that a Schur reflect is a WIP iff it is indexed
by a hook. On the other hand every WIP can be written as an expression
in the Schur reflect basis. It turns out that these expressions are
both remarkable and simple and involve only hook Schur reflects.

\begin{thm}\label{basis}
\begin{equation}\label{pins}
P_{n,\omega}=\sum_{i=0}^{n-1}(-1)^{i+1}(\omega_i-\omega_{i+1})\hat{S}_{(n-i,1^i)}, \qquad \hbox{where } \omega_0=0.
\end{equation}
\end{thm}

\noindent
{\bf Proof} 
Since each $\hat{S}_{(n-i,1^i)}$ is a weighted polynomial of weight
$(-1)^i(0,0..1,1,..)$ with the first $1$ in the $(i+1)$-th position, we
obtain that the right hand side is a weighted polynomial of level $n$
with weight vector:
$$
 \sum_i (-1)^{i+i+1}(\omega_i-\omega_{i+1})(0,0,...1,1,...)=
 (\omega_1,\omega_2,...)=\omega.
$$
\hfill
$\Box$

\section{Recursion Properties, generating functions and bases}

In~\cite{TM3} it was shown that the CSP reflects and the PSP reflects 
satisfy
Newton identities, that is, are recursive. This is in fact a property
possessed by all WIPs.

\begin{thm}\label{rec}
Let $P_{n,\omega}=P_{k,n,\omega} \in \mathfrak{F}_{k,\omega}$, then
\begin{equation}\label{recformula}
P_{n,\omega}=t_1P_{n-1,\omega}+\cdots +t_{n-1}P_{1,\omega}+t_n\omega_n
\end{equation}
\end{thm}

\noindent
{\bf Proof}
This is, essentially, the lattice definition that assigns coefficients
to monomials in a weighted family.
 Let $A(\alpha)t^\alpha$ be a monomial in $P_{n,\omega}$ and consider
 all monomials $t^\beta$ in $P_{j,\omega}$'s on the right hand side such that
 $t_{n-j}t^\beta =t^\alpha$. In the lattice $\mathcal{L}(t^\beta)$ we
 need only consider the nodes that contain $t_j$ such that $j=1,\ldots
 k$ and only those nodes $\partial_jt^\alpha$ for which $
 t_i\partial_j t^\alpha=t^\alpha$. But these are just the nodes of depth $(\sum
 \alpha_i -1)$, which by Theorem~\ref{coeffthm}, are those whose
 coefficient sum is the coefficient of $t^\alpha$.

On the other hand, every vector $\alpha$ that occurs in a monomial
$P_{j,\omega}$ in $\sum t_{n-j}P_{j,\omega}$ occurs as a vector for some
monomial in  $P_{n,\omega}$.
\hfill
$\Box$
\vskip 12pt
\begin{thm}\label{genfct}
A generating function for the WIPs in the family $\mathfrak{F}_\omega$
is
$$
\Omega(y)=\frac{\omega_1t_1y+\omega_2t_2y^2+\omega_3t_3y^3+\ldots}{1-p(y)},
\quad \mbox{where}\ p(y)= t_1y+t_2y^2+t_3y^3+\ldots
$$
\end{thm}

The polynomial $f(1/y)=x^k-t_1x^{k-1}-\cdots -t_k$, with $x=1/y$ will
be called the {\it core polynomial}. The significance of the core
polynomial will be discussed below and in the next section.
\vskip12pt
We shall organize the proof of Theorem~\ref{genfct} around the
following two lemmas.

\begin{lemma}\label{gf1}
A generating function for the family $\mathfrak{F}_{\omega_{(0)}}$
where $\omega_{(0)}=(1,1,\ldots,1,\ldots)$ is the function
$H(y)=\displaystyle{ \frac{1}{1-p(y)} }$, where $p(y)=
(t_1y+t_2y^2+t_3y^3+\ldots)$.
\end{lemma}

\noindent
{\bf Proof}
This is a consequence of Theorem 3.1 in~\cite{TM1} where it was shown that
$\{ F_n(t) \}$ is the ``positive'' multiplicative function in the core
group (under the standard convolution product in the ring of
arithmetic functions\footnote[6]{ If $\phi$ and $\theta$ are two
  arithmetic functions their convolution product is defined by $\phi *
  \theta (n)= \sum_{d|n} \phi(d)\theta(n/d) $.} ) of the group of units of the
ring of arithmetic functions, which is determined by the polynomial
$$
x^k-t_1x^{k-1}-\cdots -t_k,
$$
\noindent
which itself is the generating function for the ``negative'' sequence
with respect to the positive sequence . Thus letting $x=1/y$ gives the result.
\hfill
$\Box$

\begin{lemma}\label{gf2}
Let $\omega$ be an arbitrary weight vector, then
$P_{n,\omega}=\omega_nt_n*F_n$, where $\omega_nt_n*F_n=
\sum_{j=0}^n\omega_jt_jF_{n-j}$ and   $F_n \in
\mathfrak{F}_{\omega_{(0)}}$.
\end{lemma}

(The $*$-product is discussed in the next section where it is called
the ``level product'').
\noindent
{\bf Proof} 
An easy induction using Theorem~\ref{rec} gives the result.
\hfill
$\Box$
\vskip12pt
\noindent
{\bf Proof of Theorem~\ref{genfct} } 
 It is clear that $\sum \omega_nt_ny^n$ is the generating function for
 $\omega_nt_n$. Using Lemma~\ref{gf1} and Lemma~\ref{gf2}, we obtain
 the result we want by multiplying the generating functions.
\hfill
$\Box$

\section{Root polynomials and convolutions}

This section is motivated by work which appeared in~\cite{CG} and
~\cite{TM3} and~\cite{V}. In ~\cite{TM1} the subgroup generated by the
completely multiplicative arithmetic functions in the group of units
of the ring of arithmetic functions (the core subgroup), where
multiplication is convolution, is discussed. Each multiplicative
function in this group is uniquely determined locally (at primes) by a
particular polynomial $f(x, \boldsymbol{a})=
x^k-a_1x^{k-1}-...-a_k$. This is the {\it core polynomial} (see
Section 3, particularly Theorem~\ref{genfct}), where the parameters are 
evaluated
at $(a_1,a_2,...a_k)$. What was referred to as the ``negative'' of the
multiplicative function in that paper is an arithmetic function whose
values are just the coefficients of the core polynomial. It was proved
in ~\cite{TM1} that the ``positive'' part was determined locally by
the values of the isobaric reflects of the CSPs, that is by the GFPs
at these coefficients, or rather by ``truncations'' of the GFPs, that
is, GFPs parametrized by partitions of $n$ into parts no one of which
is greater than a fixed $k$. (Truncation is equivalent to setting the
isobaric generators $t_{k+1}, t_{k+2},...$ equal to zero in the GFPs,
or in any isobaric family).

Carrol and Gioia~\cite{CG} gave a numerical description of the rational 
roots of
the core functions. In this part of the paper, we consider isobaric
polynomials with rational coefficients and find among them polynomials
which play the same role for the rational roots of the multiplicative
arithmetic functions in the core group as the GFPs play for the core
group itself, embedding the core group into a divisible group. We do
more. Given a family of WIPs, we shall provide each of its members
with unique rational roots induced by convolution. We shall call the
products that produce these roots, {\it level products}, since they
are products of polynomials of the same level which preserve level.

The property of GFPs described above with respect to the core group in
the convolution ring of arithmetic functions can be stated as
follows. Let $p$ be a fixed prime, and let $\chi$ be a positive
element in the core, then $\chi(p^n)=F_{k,n}(a_1,...a_k)$ where
$1,a_1,...a_k$ are the coefficients of the $k$-th degree determining
polynomial. We can think of the determining polynomial in the case of   
$F_{k,n}(t_1,...t_k)$ as the generic $k$-th degree polynomial
$x^k-t_1x^{k-1}-...-t_k$. Recall that the convolution product of two
multiplicative function is given locally by
$\displaystyle{\chi_1*\chi_2(p^n)=\sum_{i=0}^n
  \chi_1(p^i)\chi_2(p^{n-i}) }$.

By induction, the $s$-th convolution power of a multiplicative
arithmetic function $\chi$ is given by $\displaystyle{ \chi^{*s}=
  \sum_{\alpha \vdash_k n}C_s(\alpha)\prod_{i=1}^k \chi(p^i)^{\alpha_i}
  }$, where $\displaystyle{ C_s(\alpha)= { {s} \choose
    {\alpha_1,...\alpha_k, (s-\sum \alpha_i)} }  }$.

We extend the definition of the convolution product for two sequence
of polynomials $\{ P_n\}_{n \geq 0}$ and $\{ Q_n\}_{n \geq 0}$ yielding
another sequence $\{ R_n\}_{n \geq 0}$ with
\begin{equation}\label{conv}
R_n:= \sum_{i=0}^n P_iQ_{n-i}.
\end{equation}
We will call the convolution product of two sequences of isobaric 
polynomials a {\it level product}, since the level is preserved.
Let $P_n= \sum_\alpha A(\alpha)t^\alpha$ be an isobaric polynomial. If $P_n$ belongs to a
weighted family the coefficients $A(\alpha)$ are given by 
Theorem~\ref{coeffthm}.

Let $q \in \mathbb{Q}$ (the group of rationals). Define the sequences 
$B_{j}^q=q(q+1)...(q+j)$ and $B_{-(j)}^q=q(q-1)...(q-j)$, for $j \geq 0$
otherwise both $B_{j}^q,\  B_{-(j)}^q$ are  zero.

\begin{thm}\label{rootf}
Let $H_{k,n}(t,q)$ denote the $q$-th convolution root of $F_{k,n}(t)$,
where $F_{k,n}(t)\in GFP$, then
\begin{equation}
 H_{k,n}(t,q)=\sum_{\alpha \vdash n} \frac{1}{(\alpha_1+\ldots
 +\alpha_k)!} B_{(\sum\alpha_i-1)}^q{ {\sum \alpha_i} \choose {\alpha_1,...\alpha_k}}t^\alpha.
\end{equation}
\end{thm}
\vskip12pt

\begin{cor}
$$
H_{k,n}(t,1)=F_{k,n}(t)
$$
\end{cor}

\noindent
{\bf Proof of Corollary}
When $q=1$, then $B_{(\sum\alpha_i-1)}^1=(\sum_{i=1}^k \alpha_i)!$. The
Corollary then follows from Theorem~\ref{coeffthm} (see remark
following that theorem).
\hfill
$\Box$
\vskip12pt

\begin{cor}
If $\chi$ belongs to the Core of the group of units of the convolution
ring of arithmetic functions, then $H_{k,n}(\boldsymbol{a},q)=
\chi^{*q}(p^n)$, where $\boldsymbol{a}= (a_1,\ldots a_k)$ is the set of
coefficients of the core polynomial of $\chi$.
\end{cor}
\hfill
$\Box$

The proof of Theorem~\ref{rootf} follows from a more general
result. Namely that each polynomial in a weighted family of isobaric
polynomials has a unique $q$-root for every rational number $q$.

\begin{thm}\label{hgeneral}
$\displaystyle{ H_{k,n,\omega}= \sum_{\alpha \vdash_k
    n}L_{k,n,\omega}(\alpha)t^\alpha }$ where
$$
 L_{k,n,\omega}(\alpha)= \sum_{j=0}^{\sum \alpha_i-1} \frac{1}{(\prod
   \alpha_i)!}
 {{\sum \alpha_i-1} \choose j} B_{-(j)}^qD_{(\sum\alpha_i-j-1)}(\omega_1^{\alpha_1}\ldots \omega_k^{\alpha_k})
$$
is the $q$-th level root of $P_{k,n,\omega} \in
\mathfrak{F}_{k,\omega}$ and $H_{k,0,\omega}=1$.
\end{thm}

Theorem~\ref{rootf} follows from Theorem~\ref{hgeneral}. The following
lemmas will be used in proving this.

\begin{lemma}
$$
D_j(\omega^\alpha)|_{\omega=(1,1,..)}= \frac{(\sum \alpha_i)!}{(\sum\alpha_i-j)!}.
$$
\end{lemma}

\noindent
{\bf Proof}
At depth $0$ the value is $1$. If after $(j-1)$ derivations the value
is $\displaystyle{\frac{(\sum \alpha_i)!}{(\sum\alpha_i-j+1)!} }$,
then in the $j$-th step the exponent sum is decreased by $1$, so by derivation,
$(\sum\alpha_i-j)$ appears as the only new factor in the value of
$D_j\omega^\alpha$.
\hfill
$\Box$ 

\begin{lemma}\label{stir}
$$
\sum_{j=0}^{\sum \alpha_i-1} \frac{(\sum\alpha_i-1)!}{(j+1)!} { {\sum \alpha_i -1} \choose j }B_{-(j)}^q= B^q_{(\sum \alpha_i-1)}.
$$
\end{lemma}

\noindent
{\bf Proof}
Consider the following Stirling functions $[x]_n=x(x-1)\ldots
(x-n+1)$, and $[x]^n=x(x+1)\ldots (x+n-1)$. From the theory of
Stirling numbers of 1st and 2nd kind we have the relation
$\displaystyle{ [x]^p= \sum_{j=1}^{p}
  {{p-1} \choose {j-1}} \frac{p!}{j!}[x]_{j} }$ (e.g.,~\cite{T}, p15, problem 3). This translates into
$$
B_{(p)}^q= \sum_{j=0}^p {p \choose j}\frac{(p+1)!}{(j+1)!}B_{-(j)}^q.
$$ 
Letting now $p=\sum \alpha_i-1$ in the Equation above gives the result
we wanted.
\hfill
$\Box$
\vskip12pt
Theorem~\ref{hgeneral} is a consequence of the following theorem,
which shows an interesting closure property.

\begin{thm}\label{hsum}
$$
H_{k,n,\omega}(t,q)*H_{k,n,\omega}(t,q')=H_{k,n,\omega}(t,q+q')
$$
\end{thm}
\vskip12pt
Before we proceed we need these lemmas.

\begin{lemma}\label{l1}
$$
\sum_{j=0}^{n+1}{ {n+1} \choose j} B_{-(n-j)}^qB_{-(j-1)}^{q'} = B_{-(n)}^{q+q'}
$$
\end{lemma}

\noindent
{\bf Proof} As in Lemma~\ref{stir} this is a consequence of the theory of
Stirling numbers. Here the relevant result is that $[x+y]_{n+1}=
\sum_{j=0}^{n+1}{ {n+1} \choose j }[x]_{n+1-j}[y]_j$.
An analogous formula for $[x+y]^{n+1}$ shows that 
$\sum_{j=0}^{n+1}{ {n+1} \choose j }B_{(n-j)}^qB_{(j-1)}^{q'} =
B_{(n)}^{q+q'}$.
\hfill
$\Box$
\vskip12pt

\begin{lemma}\label{l2}
Let $\alpha=(\alpha_1,\ldots \alpha_k)$ with $\alpha_i \geq 0$ and
$|\alpha|=\sum_{i=1}^k \alpha_i=n$. Then
\begin{equation}\label{deriv} 
 \sum_{|\beta|= m,\,\beta \leq \alpha} \prod {\alpha_i \choose 
\beta_i}D_p\omega^\beta
D_q\omega^{\alpha-\beta} = {{n-p-q}\choose{m-p}} D_{p+q} \omega^\alpha
\end{equation}
where $m \leq n$ and $p,\, q \in \mathbb{N}$.
\end{lemma}

\noindent
{\bf Proof}
We will prove this lemma by induction.
 
Let $\mathcal{P}(q)$ be the following statement:``(\ref{deriv}) is
true for every $p$''. First we will show that $\mathcal{P}(0)$ is
true, i.e.
\begin{equation}\label{p0}
\sum_{|\beta|= m,\,\beta \leq \alpha} \prod {\alpha_i \choose
  \beta_i}(D_p\omega^\beta )\omega^{\alpha-\beta} = {{n-p}\choose{m-p}} D_{p} \omega^\alpha
\end{equation}

\noindent
Before we proceed let us note the case where $p=0$, which will be used
in the sequel. In this case the identity ~(\ref{p0}) becomes
$\sum_\beta \prod_i {\alpha_i \choose
  \beta_i}\omega^\beta\omega^{\alpha-\beta} =
{{n}\choose{m}}\omega^\alpha$, that is $\sum_\beta\prod_i {\alpha_i \choose\beta_i} ={{n}\choose{m}}$. The
latter is true since the left hand side of it is the coefficient of
$x^m$ in the product
$(1+x)^{\alpha_1}...(1+x)^{\alpha_k}=(1+x)^{\sum\alpha_i}=(1+x)^n$,
which is clearly $n \choose m$.

For the general case we may write $D_p\omega^{\beta}=\sum_{|\beta-\gamma|=p} f_\gamma(\beta) \omega^\gamma$, where
 $f_\gamma(\beta)$ is a polynomial in
 the  $\beta$'s. More precisely, if $s_i= \beta_i-\gamma_i$ then
$f_\gamma(\beta)= c\cdot\prod [\beta_i]_{s_i}=c \prod
  {{\beta_i}\choose {\gamma_i}}(\beta_i-\gamma_i)!$, where $c$ is a
  constant. This constant counts the number of paths in the lattice
  $\mathcal{L}(t^\alpha)$ from $t^\beta$ to $t^\alpha$. Equivalently,
  $c$ is the enumeration of $\{\alpha_1,...\alpha_1-s_1+1,..\}$ such
  that $\alpha_1$ precedes $\alpha_1-1$ and so on. 
Hence $c= \frac{(s_1+...+s_k)!}{s_1!...s_k!}= \frac{p!}{\prod
  (\beta_i-\gamma_i)!}$. In fact we have showed that
\begin{equation}\label{diffp}
D_p\omega^{\beta}=p!\sum_{|\beta-\gamma|=p} \prod {{\beta_i} \choose {\gamma_i}}\omega^\gamma.
\end{equation}

We can now rewrite Equation~\ref{p0} as
$$
\sum_{|\beta|= m} \prod {\alpha_i
  \choose\beta_i}p!\sum_{|\beta-\gamma|=p} \prod {{\beta_i} \choose
  {\gamma_i}}\omega^\gamma \omega^{\alpha-\beta} = {{n-p}\choose{m-p}} p!\sum_{|\alpha-\delta|=p} \prod {{\alpha_i} \choose {\delta_i}}\omega^\delta
$$
\noindent
For $\delta$ fixed on the right hand side we must have equality of the
corresponding coefficients, i.e.
$$
\sum_{|\beta|= m} \prod {{\alpha_i}\choose{\beta_i}} {{\beta_i}
  \choose {\delta_i-\alpha_i+\beta_i}}={{n-p}\choose{m-p}}\prod {{\alpha_i} \choose {\delta_i}}
$$
\noindent
and if we rewrite the left hand side we obtain
$$
\sum_{|\beta|= m} \prod {{\alpha_i} \choose {\delta_i}} {{\delta_i}
  \choose {\alpha_i-\beta_i}}={{n-p}\choose{m-p}}\prod {{\alpha_i} \choose {\delta_i}}
$$
\noindent 
which gives the identity we showed for $p=0$.

To show the induction step we differentiate the expression in
$\mathcal{P}(q)$ (~\ref{deriv}) to get
$$
\sum_{|\beta|=m} \prod_i {\alpha_i \choose
  \beta_i}D_{p+1}\omega^\beta D_q\omega^{\alpha-\beta} + \sum_{|\beta|= m} \prod_i {\alpha_i \choose
  \beta_i}D_{p}\omega^\beta D_{q+1}\omega^{\alpha-\beta} = {{n-p-q}\choose{m-p}} D_{p+q+1} \omega^\alpha
$$

and by the induction step
$$
\sum_{|\beta|=m} \prod_i {\alpha_i \choose
  \beta_i}D_{p}\omega^\beta D_{q+1}\omega^{\alpha-\beta}= \Big{[}
{{n-p-q}\choose{m-p}}-{{n-p-q-1}\choose{m-p-1}}\Big{]} D_{p+q+1}
\omega^\alpha
$$
$$
={{n-p-q-1}\choose{m-p}} D_{p+q+1} \omega^\alpha
$$
\noindent
which is exactly $\mathcal{P}(q+1)$, and thus the proof of
Lemma~\ref{l2} is complete.
\hfill
$\Box$
\vskip12pt
We need one more lemma which is a binomial identity.
\begin{lemma}\label{l3}
Let $p,\ q$ and $n$ such that $p+q <n$. Then
$$
\sum_{i=0}^{n-p-q} {{p+i} \choose p} {{n-p-i} \choose q} = {{n+1}\choose
{p+q+1}}
$$
\end{lemma}

\noindent
{\bf Proof}
If we make the convention that ${a \choose b}=0$ if $a <b$, we can
think of the left hand side as being the sum $\displaystyle{
  \sum_{i=0}^n {i \choose p}{{n-i} \choose q}}$.
This is in fact the coefficient of $x^py^q$ in 

$$
\sum_{i=0}^n (1+x)^i(1+y)^{n-i} = \frac{(1+x)^{n+1}-(1+y)^{n+1}}{x-y}.
$$
\noindent
Let us denote by $a(i,j)$ the coefficient of $x^iy^j$ in the
expression above, i.e. 
$$
\frac{(1+x)^{n+1}-(1+y)^{n+1}}{x-y}= \sum_{i,j}a(i,j)x^iy^j.
$$
We have that 
$$ 
(1+x)^{n+1}-(1+y)^{n+1}=\sum_{i,j}a(i,j)(x^{i+1}y^j -x^iy^{j+1}) 
=\sum_{i,j} [a(i-1,j)-a(i,j-1)]x^iy^j.
$$
\noindent
Here $a(i,j)=0$ if one of $i,\ j$ is negative. By equating
coefficients, we obtain
$$
 a(i,0)=a(0,i)= {{n+1} \choose {i+1}}
$$
\noindent
and if both $i,\, j >0$, then we get $a(i-1,j)=a(i,j-1)$. An easy
inductive argument shows that
$a(i,j)=a(i+1,j-1)=...=a(i+j,0)=\displaystyle{ {{n+1} \choose {i+j+1}}
  }$ and the proof of the lemma is complete.
\hfill
$\Box$
\vskip12pt  
\noindent
{\bf Proof of Theorem~\ref{hsum} }
By the definition of the level product we have that
$$
H_{n,\omega}(t,q)*H_{n,\omega}(t,q')=\sum_{i=0}^n
H_{n-i,\omega}(t,q)H_{i,\omega}(t,q')= \sum_{i=0}^n(\sum_{\beta \vdash
  n-i}L^q(\beta)t^\beta)(\sum_{\gamma \vdash i}L^{q'}(\gamma)t^\gamma)
$$
$$
= \sum_{\alpha \vdash n} \sum _{\beta \leq \alpha}
L^q(\beta)L^{q'}(\alpha-\beta)t^\alpha
$$
\noindent
where $\beta \leq \alpha$ means $\beta_i \leq \alpha_i$, for every
$i$.
Therefore we need to show that
$$
\sum _{\beta \leq \alpha} L^q(\beta)L^{q'}(\alpha-\beta) =
L^{q+q'}(\alpha).
$$
By replacing $L$'s with their formulas (see definition in the
statement of Theorem~\ref{hgeneral}),  we therefore need to show that
$$
 \sum_{\beta \leq \alpha} \prod_i {{\alpha_i} \choose {\beta_i}} \sum_{s=0}^{r-1}\sum_{t=0}^{p-r-1} { {r-1} \choose s}
{ {p-r-1} \choose t} B_{-(s)}^qB_{-(t)}^{q'}
D_{(r-s-1)}\omega^\beta
D_{(p-r-t-1)}\omega^{\alpha-\beta}
$$
\begin{equation}\label{sides}
= \sum_{j=0}^{p-1}{{p-1} \choose j}
B_{-(j)}^{q+q'} D_{(p-j-1)}\omega^\alpha
\end{equation}
\noindent
where we denote for simplicity $|\alpha|=p$ and $|\beta|=r$.
Fix now an index $j$ in the right hand side above. An important
fact is that for each such $j$ the expression of
$D_{p-j-1}\omega^\alpha$ gives an homogeneous polynomial of degree
$(j+1)$ in $\omega_1,...\omega_k$. So it suffices to show that the
two corresponding homogeneous polynomials of the same degree on both sides
coincide. To obtain the homogeneous polynomial of degree $(j+1)$ on
the left hand side we need to pick indices $s$ and $t$ such that
$(s+1)+(t+1)=j+1$ i.e. $t=j-s-1$. The homogeneous polynomial of degree
$(j+1)$ on the left hand side is therefore
$$
 \sum_{\beta \leq \alpha} \prod_i{{\alpha_i} \choose {\beta_i}} \sum_{s=0}^{r-1}{ {r-1} \choose s}
{ {p-r-1} \choose {j-s-1} } B_{-(s)}^qB_{-(j-s-1)}^{q'}
D_{(r-s-1)}\omega^\beta
D_{(p-r-j+s)}\omega^{\alpha-\beta}
$$    
As before we consider ${a \choose b}=0$ if $a<b$.
We can rewrite $\sum_{\beta \leq \alpha}$ as $\sum_{r=0}^p
\sum_{|\beta|= r,\, \beta \leq \alpha}$, and the expression above
becomes 
$$
\sum_{r=0}^p\sum_{|\beta|=r,\,\beta \leq \alpha} \prod_i {
  {\alpha_i} \choose {\beta_i}} \sum_{s=0}^{r-1} { {r-1} \choose s}
{ {p-r-1} \choose {j-s-1}} B_{-(s)}^qB_{-(j-s-1)}^{q'}
D_{(r-s-1)}\omega^\beta
D_{(p-r-j+s)}\omega^{\alpha-\beta}
$$
$$
=\sum_{r=0}^p\sum_{s=0}^{r-1} { {r-1} \choose s}
{ {p-r-1} \choose {j-s-1}} B_{-(s)}^qB_{-(j-s-1)}^{q'}
\sum_{|\beta|=r,\,\beta \leq \alpha} \prod_i {
  {\alpha_i} \choose {\beta_i}} D_{(r-s-1)}\omega^\beta
D_{(p-r-j+s)}\omega^{\alpha-\beta}
$$
\noindent 
which by using Lemma~\ref{l2} is
$$
= \sum_{r=0}^p\sum_{s=0}^{r-1} { {r-1} \choose s}
{ {p-r-1} \choose {j-s-1}} B_{-(s)}^qB_{-(j-s-1)}^{q'} { {j+1}\choose
  {s+1}} D_{p-j-1}\omega^\alpha
$$
$$
= \sum_{s=0}^{p-1}{ {j+1}\choose {s+1}}  B_{-(s)}^qB_{-(j-s-1)}^{q'}
D_{p-j-1}\omega^\alpha \sum_{r=0}^p{ {r-1} \choose s}
{ {p-r-1} \choose {j-s-1}} 
$$
and by Lemma~\ref{l3} is
$$
=\sum_{s=0}^{p-1}{ {j+1}\choose {s+1}}  B_{-(s)}^qB_{-(j-s-1)}^{q'}
D_{p-j-1}\omega^\alpha {{p-1} \choose j} 
$$
\noindent
which finally by Lemma~\ref{l1} is
$$
={{p-1} \choose j} B_{-(j)}^{q+q'} D_{p-j-1}\omega^\alpha 
$$    
and this is exactly the homogeneous polynomial of degree $(j+1)$ in the
right hand side of~(\ref{sides}) and thus the proof of 
Theorem~\ref{hsum} is complete.

\hfill
$\Box$
\section{Algebraic structure}

Let $\mathcal{H}_{n,\omega}$ denote the algebra generated by all $H_{n,\omega}$
under the addition and the level product. As a consequence of
Theorem~\ref{hsum} each $H_{n,\omega}$ has a level-product inverse in $\mathcal{H}$.

\begin{thm}\label{inverse}
$$
H_{n,\omega}^{-1} (t,q) = H_{n,\omega}(t, -q)
$$
\end{thm}
\hfill
$\Box$
\vskip12pt

So from Theorems~\ref{hsum} and ~\ref{inverse} we have that for a fixed 
weight $\omega$ and a given level $n$ the polynomials  
$\mathcal{H}_{n,\omega}= \{ H_{n,\omega}(t,q) \}_{q \in \mathbb{Q}}$ form 
an abelian group under the level product isomorphic to the rationals, $\mathbb{Q}$,
under addition. The group $\mathcal{T}= \{ \mathfrak{F}_{\omega} \}_{\omega}$ 
acts on this group by  translation in the following way: 
$\mathcal{T}$ acts on a family of WIPs by (say) a right translation 
(Theorem~\ref{free}) and in a natural way the $q$-th roots follow along. 
Theorem~\ref{inverse} applies to a family of WIPs, as well, giving the 
subgroup determined by a weighted isobaric family under the level 
operation. All of this together with the derivation operators $\partial_j$ 
give a structure of differential graded group to $\mathcal{H}= 
\oplus_\omega\oplus_n \mathcal{H}_{n,\omega}$ acted on by an affine group. 

\section{Appendix}

{\it Schur reflects for $Sym(n)$, with $n=1,2,...6$. }

\vskip12pt

$
\begin{array}{llr}
\hat{S}_{(1)} & = & t_1
\end{array}
$
\vskip12pt
$
\begin{array}{llr}
\hat{S}_{(2)} & = & t_1^2+t_2\\
\hat{S}_{(1^2)} & = & -t_2 
\end{array}
$
\vskip 12pt
$
\begin{array}{llr}
\hat{S}_{(3)} & = & t_1^3+2t_1t_2+t_3\\
\hat{S}_{(2,1)} & = & -t_1t_2 -t_3\\
\hat{S}_{(1^3)} & = & t_3 
\end{array}
$
\vskip12pt
$
\begin{array}{llr}
\hat{S}_{(4)} & = & t_1^4+3t_1^2t_2+t_2^2+ 2t_1t_3+ t_4\\
\hat{S}_{(3,1)} & = & -t_1^2t_2-t_2^2-t_1t_3 -t_4\\
\hat{S}_{(2^2)} & = & t_2^2-t_1t_3\\
\hat{S}_{(2,1^2)} & = & t_1t_3 +t_4\\
\hat{S}_{(1^4)} & = & -t_4
\end{array}
$
\vskip12pt
$
\begin{array}{llr}
\hat{S}_{(5)} & = & t_1^5+4t_1^3t_2+3t_1t_2^2+ 3t_1^2t_3+2t_2t_3+2t_1 t_4+t_5\\
\hat{S}_{(4,1)} & = & -t_1^3t_2-2t_1t_2^2-t_1^2t_3 -2t_2t_3
-t_1t_4-t_5\\
\hat{S}_{(3,2)} & = & t_1t_2^2-t_1^2t_3 +t_2t_3-t_1t_4\\
\hat{S}_{(3,1^2)} & = &t_1^2t_3+t_2t_3+t_1t_4+t_5\\ 
\hat{S}_{(2^2,1)} & = & -t_2t_3+t_1t_4\\
\hat{S}_{(2,1^2)} & = & t_1t_3 +t_4\\
\hat{S}_{(2,1^3)} & = & -t_1t_4 -t_5\\
\hat{S}_{(1^5)} & = & t_5
\end{array}
$
\vfill

\pagebreak
$
\begin{array}{llr}
\hat{S}_{(6)} & = & t_1^6+5t_1^4t_2+6t_1^2t_2^2+ t_2^3 +4t_1^3t_3
+t_3^2 + 6t_1t_2t_3+3t_1^2t_4+2t_2t_4+2t_1t_5+  t_6\\
\hat{S}_{(5,1)} & = & -t_1^4t_2- 3t_1^2t_2^2-t_2^3-t_1^3t_3-4t_1t_2t_3-t_3^2-2t_2t_4 -2t_1t_5-t_6\\
\hat{S}_{(4,2)} & = & t_1^2t_2^2+t_2^3 -t_1^3t_3-t_1^2t_4 +t_2t_4-t_1t_5\\
\hat{S}_{(4,1^2)} & = &
t_1^3t_3+2t_1t_2t_3+t_3^2+t_1^2t_4+t_2t_4+t_1t_5+t_6\\
\hat{S}_{(3^2)} & = & 2t_1t_2t_3-t_2^3+t_3^2-t_1^2t_4-t_2t_4\\ 
\hat{S}_{(3,2,1)} & = & -t_1t_2t_3-t_3^2 +t_1^2t_4+t_1t_5 \\
\hat{S}_{(2^3)} & = & t_3^2-t_2t_4\\
 \hat{S}_{(3,1^3)} & = & -t_1^2t_4-t_2t_4 -t_1t_5-t_6\\
\hat{S}_{(2^2,1^2)} & = & t_2t_4- t_1t_5\\
\hat{S}_{(2,1^4)} & = & t_1t_5 +t_6\\
\hat{S}_{(1^6)} & = &- t_6
\end{array}
$

\end{document}